\documentclass{amsart}
\usepackage{graphics}
\usepackage{amssymb}

\newcommand{\floor}[1]{\left\lfloor{#1}\right\rfloor}

\newcommand{\abs}[1]{\left\lvert{#1}\right\rvert}
\newcommand{\norm}[1]{\left\|{#1}\right\|}
\newcommand{\slo}{\operatorname{slo}}

\newcommand{\bb}{\mathbb}
\newcommand{\ol}{\overline}

\newcommand{\pr}{\mathrm{pr}}
\newcommand{\GL}{\mathrm{GL}}
\newcommand{\SL}{\mathrm{SL}}
\newcommand{\R}{\mathbb{R}}
\newcommand{\N}{\mathbb{N}}
\newcommand{\Z}{\mathbb{Z}}
\newcommand{\Q}{\mathbb{Q}}
\newcommand{\T}{\mathbb{T}}

\newcommand{\ZCP}{\mathbb{Z}^2_{\mathrm{cp}}}

\DeclareMathOperator{\rot}{\rho}
\DeclareMathOperator{\cl}{cl}
\DeclareMathOperator{\diam}{diam}

\newcommand{\homeo}{\operatorname{Homeo}}
\newcommand{\homi}{{\operatorname{Homeo}}_*}

\newcommand{\cf}{c.f.\ }
\newcommand{\ie}{i.e.\ }

\newcommand{\proc}[1]{\medbreak\noindent{\it #1}\hspace{1ex}\ignorespaces}
\newcommand{\ep}{\qed}
\newtheorem{theorem}{Theorem}[section]

\newtheorem{lemma}[theorem]{Lemma}
\newtheorem{proposition}[theorem]{Proposition}
\newtheorem{conjecture}[theorem]{Conjecture}

\newcommand{\acks}{\subsection*{Acknowledgements}}
\newcounter{thmA}

\newtheorem{theoremA}[thmA]{Theorem}

\newcounter{thmB}

\newtheorem{theoremB}[thmB]{Theorem}

\newcounter{thmN}

\newtheorem{theorem*}[thmN]{Theorem}

\newcommand{\definition}{\addtocounter{theorem}{1}\proc{Definition \thetheorem.}}
\newcommand{\definitionend}{\medbreak}
\newcommand{\remark}{\refstepcounter{theorem}\proc{Remark \thetheorem.}}
\newcommand{\remarkend}{\medbreak}

\begin{document}

\title{Free curves and periodic points for torus homeomorphisms}
\author[A. Kocsard]{Alejandro Kocsard}

\address{Instituto Nacional de Matem\'atica Pura e Aplicada, Estrada
  Dona Castorina, 110 Rio de Janeiro, Brasil}

\email{alejo@impa.br}

\thanks{The first author was supported by FAPERJ-Brasil.}

\author[A.  Koropecki]{Andr\'es Koropecki}

\address{Universidade Federal Fluminense, Instituto de Matem\'atica, Rua M\'ario Santos Braga S/N, 24020-140 Niter\'oi, RJ, Brasil}

\email{akoro@impa.br}

\thanks{The second author was supported by CNPq-Brasil.}

\begin{abstract}
  We study the relationship between free curves and periodic points
  for torus homeomorphisms in the homotopy class of the identity. By free curve we mean a homotopically nontrivial simple closed curve that is disjoint from its image.
We prove that every rational point in the rotation set is realized by a
  periodic point provided that there is no free curve and the rotation
  set has empty interior.  This gives a topological version of a
  theorem of Franks. Using this result, and inspired by a theorem of
  Guillou, we prove a version of the Poincar\'e-Birkhoff Theorem for
  torus homeomorphisms: in the absence of free curves, either there is
  a fixed point or the rotation set has nonempty interior.
\end{abstract}
\maketitle

\section{Introduction}

Given a torus homeomorphism $F\colon\bb{T}^d\to\bb{T}^d$ in the
identity homotopy class, the \emph{rotation set} of a lift
$f\colon\bb{R}^d\to\bb{R}^d$ of $F$ was introduced by Misiurewicz and
Ziemian in \cite{m-z}, and it is defined as the set of all accumulation
points of sequences of the form
\begin{equation*}
  \left\{\frac{f^{n_i}(x_i)-x_i}{n_i}\right\}_{i\in\N}
\end{equation*}
where $n_i\to \infty$ and $x_i\in \R^2$.

This set carries dynamical information about $F$, but understanding
this information is not as easy as in the one-dimensional case.
Moreover, when $d=2$ the rotation set has nice geometric properties
(for example, convexity) which are no longer valid in
higher dimension. For this reason, we restrict our attention to the
two-dimensional setting.

A central problem, inspired by the Poincar\'e theory for circle
homeomorphisms, is to determine when a point $(p_1/q,p_2/q)\in
\rot(f)\cap \Q^2$ (with $\gcd(p_1,p_2,q)=1$) is realized by a periodic
orbit of $F$, \ie when can we find a point $x\in \R^2$ such that
\begin{equation*}
  f^q(x) = x + (p_1,p_2).
\end{equation*}

This problem has been thoroughly studied, especially by Franks, who
proved in \cite{franks:reali-ext} that extremal points of the rotation
set are always realized by periodic orbits. However, it is generally
not true that \emph{every} rational point in the rotation set is
realized.  In this aspect, the case that is best understood is when
the rotation set has non-empty interior. In fact, Franks proved in
\cite{franks:reali-int} that every rational point in the interior of
the rotation set is realized by a periodic orbit, and this is optimal
in the sense that there are examples (even area-preserving ones) where
the rotation set has non-empty interior and many rational points on
the boundary, but the only ones that are realized by periodic orbits
are extremal or interior ones (see \cite[\textsection 3]{m-z2}).

On the other hand, when the rotation set has empty interior the
situation is more delicate. It is easy to construct a
periodic-point-free homeomorphism $F\colon \T^2\to\T^2$ such that its
rotation set is a segment containing many rational points. However,
this cannot happen under some additional hypothesis:

\begin{theorem}[Franks \cite{franks:reali-area}]
  \label{th:franks-area}
  If an area-preserving homeomorphism of $\T^2$ in the homotopy class
  of the identity has a rotation set with empty interior, then every
  rational point in its rotation set is realized by a periodic orbit.
\end{theorem}

In the same article, Franks asks whether the area-preserving
hypothesis is really necessary for the conclusion of the theorem. It
is natural to expect that a weaker, more topological hypothesis should
suffice to obtain the same result. This topological substitute for the
area-preserving hypothesis turns out to be, to some extent, the
\emph{curve intersection property}. An essential simple closed curve
is \emph{free} for $F$ if $F(\gamma)\cap \gamma=\emptyset$. We say
that $F$ has the curve intersection property if $F$ has no free
curves. Our first result can be stated as follows:
\medbreak

\begin{theoremA}
\it{  If a homeomorphism of $\T^2$ in the homotopy class of the identity
  satisfies the curve intersection property and its rotation set has
  empty interior, then every rational point in its rotation set is
  realized by a periodic orbit.}
\end{theoremA}

Our proof is essentially different of that of Franks' theorem, since
the latter relies strongly on chain-recurrence properties that are
guaranteed by the area preserving hypothesis but not by the curve
intersection property.

Variations of the curve intersection property are already present in
some fixed point theorems. An interesting case is a generalization of
the classic theorem of Poincar\'e-Birkhoff \cite{birkhoff:PB}, which
states that for an area-preserving homeomorphism of the closed
annulus, isotopic to the identity and verifying the \emph{boundary
twist condition}, there exist at least two fixed points.  Birkhoff and
Ker\'ekj\'art\'o already noted that, for getting at least one fixed
point, the area preserving hypothesis was not completely necessary,
and they replaced it by the weaker condition that any essential simple
closed curve intersects its image by the homeomorphism
\cite{kerekjarto:PB}. 

An even more topological version of this theorem was proved by
Guillou, who substituted the twist condition by the property that
every simple arc joining the boundary components intersects its image
by $F$:
\begin{theorem*}[Guillou, \cite{guillou:PB}]
  If $F\colon [0,1]\times S^1\to[0,1]\times S^1$ is an orientation-preserving
  homeomorphism isotopic to the identity and such that every essential
  simple closed curve or simple arc joining boundary components
  intersects its image by $F$, then $F$ has a fixed point.
\end{theorem*}

The hypotheses of the above theorem can be regarded as the curve
intersection property in the setting of the closed annulus. This led
Guillou to ask if a similar result holds for the torus. The answer is
no, as an example by Bestvina and Handel shows \cite{bestvina-handel}.
Their example relies in the fact that the existence of a free curve
imposes a restriction on the ``size'' of the rotation set (see Lemma
\ref{lem:CIP-3}). For this reason, their example has a rotation set
with nonempty interior, which implies that it has infinitely many
periodic points of arbitrarily high periods \cite{franks:reali-int},
and positive topological entropy \cite{llibre-mackay}. The question
that naturally arises is whether the presence of this kind of ``rich''
dynamics is the only new obstruction to the existence of free curves
in $\T^2$.  In other words, \emph {is the answer to Guillou's question
  affirmative if the rotation set has empty interior?} This leads to
our next result.

\begin{theoremB}
  Let $F\colon \T^2\to \T^2$ be a homeomorphism homotopic to the
  identity and satisfying the curve intersection property. Then either
  $F$ has a fixed point, or its rotation set has nonempty interior.
\end{theoremB}

\begin{figure}[ht]
\centering{\resizebox{0.6\textwidth}{!}{\includegraphics{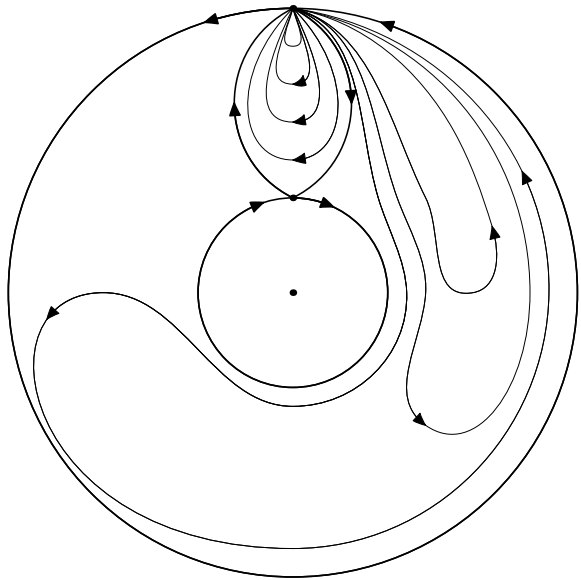}}}
  \caption{A curve-intersecting flow with no periodic points of period $> 1$} \label{fig:twofixed}
\end{figure}

Thus, if $F$ has the curve intersection property, then either $F$ has
periodic orbits of arbitrarily high periods or it has a fixed point.
In the latter case, one might expect the existence of a second fixed
point (as in the case of the annulus) but no more than that.  Figure
\ref{fig:twofixed} shows that one cannot expect more than two fixed
points in the annulus; the time-one map of the flow sketched there is
a homeomorphism with the curve intersection property, which has two
fixed points and no other periodic points. Gluing a symmetric copy of
this homeomorphism through the boundaries of the annulus, one obtains
an example with the same properties in $\T^2$.

\begin{conjecture} 
  Under the hypotheses of Theorem B, either $F$ has two fixed points
  or $F$ has periodic points of arbitrarily high periods.
\end{conjecture}

\section{Notation and preliminaries}
\label{sec:prelim}

As usual, we denote the $2$-torus $\R^2/\Z^2$ by $\T^2$, being
$\pi\colon \R^2\to\T^2$ the canonical quotient projection. We denote
the integer translations by
$$T_1\colon(x,y)\mapsto (x+1,y)\,\text{ and }\,
T_2\colon(x,y)\mapsto(x,y+1),$$ and $\pr_i\colon \R^2\to\R$, for
$i=1,2$ are the projections onto the first and second coordinate,
respectively.

By $\homeo(X)$ we mean the set of homeomorphisms of $X$ onto itself,
and by $\homi(X)$ the set of elements of $\homeo(X)$ which are
homotopic to the identity.

Given $F\in\homeo(\T^2)$, a lift of $F$ to $\R^2$ is a map $f\in
\homeo(\R^2)$ such that $\pi f = F\pi$. A homeomorphism $f\colon
\R^2\to\R^2$ is a lift of an element of $\homi(\T^2)$ if and only if
$f$ commutes with $T_1$ and $T_2$. Any two lifts of a given
homeomorphism of $\T^2$ always differ by an integer translation.  We
will usually denote maps of $\T^2$ to itself by uppercase letters, and
their lifts to $\R^2$ by their corresponding lowercase letters.

By $\ZCP$ we denote the set of pairs of integers $(m,n)$ such that $m$
and $n$ are coprime.  We will say that $(x,y)\in \R^2$ is an integer
point if both $x$ and $y$ are integer, and a rational point if both
$x$ and $y$ are rational. When we write a rational number as $p/q$, we
assume that it is in reduced form, \ie that $p$ and $q$ are coprime,
except when we are talking about a rational point $(p_1/q,p_2/q)$, in
which case we assume that $p_1$, $p_2$ and $q$ are mutually coprime
(\ie $\gcd\{p_1,p_2,q\}=1$).

\subsection{The rotation set.}
From now on, we will assume that $F\in \homi(\T^2)$, and $f\colon
\R^2\to\R^2$ is a lift of $F$.

\definition (Misiurewicz \& Ziemian \cite{m-z})
  The \emph{rotation set of $f$} is defined as
  \begin{align*}
    \rot(f)&=\bigcap_{m=1}^{\infty}\cl\left(\bigcup_{n=m}^\infty
      \left\{\frac{f^n(x)-x}{n}:x\in\R^2\right\}\right)\subset\R^2
  \end{align*}
  The \emph{rotation set of a point} $x\in \R^2$ is defined by
  $$\rot(f,x)= \bigcap_{m=1}^{\infty}\cl \left\{\frac{f^n(x)-x}{n}:
    n>m\right\}.$$

  If the above set consists of a single point $v$, we call $v$ the
  rotation vector of $x$.
\definitionend

\remark
  \label{pro:rot}
  It is easy to see that for integers $n,m_1,m_2$,
$$\rot(T_1^{m_1}T_2^{m_2}f^n) = n\rot(f)+(m_1,m_2).$$
In particular, the rotation set of any other lift of $F$ is an integer
translate of $\rot(f)$, and we can talk about the ``rotation set of
$F$'' if we keep in mind that it is defined modulo $\Z^2$.
\remarkend

\begin{theorem}[\cite{m-z}]
  \label{th:m-z-rot}
  The rotation set is compact and convex, and every extremal point of
  $\rot(f)$ is the rotation vector of some point.
\end{theorem}

Given $A\in \GL(2,\Z)$, we denote by $\tilde A$ the homeomorphism of
$\T^2$ lifted by it. If $H\in\homeo(\T^2)$, there is a unique
$A\in\GL(2,\Z)$ such that for every lift $h$ of $H$, the map $h-A$ is
bounded (in fact, $\Z^2$-periodic). From this it follows that $H$ is
isotopic to $\tilde A$, and $h^{-1}-A^{-1}$ is bounded. In fact,
$H\mapsto A$ induces an isomorphism between the isotopy group of
$\T^2$ and $\GL(2,\Z)$.

\begin{lemma}
  \label{lem:conjuga}
  Let $F\in\homi(\T^2)$, $A\in\GL(2,\Z)$, and $H\in\homeo(\T^2)$
  isotopic to $\tilde{A}$. Let $f$ and $h$ be the respective lifts of
  $F$ and $H$ to $\R^2$. Then
  \begin{equation*}
    \rot(hfh^{-1})=A\rot(f).
  \end{equation*}
  In particular, $\rot(AfA^{-1})=A\rot(f)$.
\end{lemma}

\proc{Proof.}
  We can write $((hfh^{-1})^n(x) - x)/n$ as
  $$\frac{(h-A)(f^nh^{-1}(x))}{n} +
  A\left(\frac{f^n(h^{-1}(x))-h^{-1}(x)}{n}\right) +
  A\left(\frac{(h^{-1}-A^{-1})(x)}{n}\right)$$ and using the fact that
  $h-A$ and $h^{-1}-A^{-1}$ are bounded, we see that the leftmost and
  rightmost terms of the above expression vanish when $n\to\infty$.
  Thus if $n_k\to\infty$ and $x_k\in\R^2$, we have
  $$\lim_{k\to\infty}\frac{(hfh^{-1})^{n_k}(x_k) - x_k}{n_k} =
  A\left(\lim_{k\to\infty}
    \frac{f^{n_k}(h^{-1}(x_k))-h^{-1}(x_k)}{n_k}\right)$$ whenever the
  limits exist. Since $h$ is a homeomorphism, it follows from the
  definition that $\rot(hfh^{-1})=A\rot(f)$.
\ep\medbreak

We will use the above lemma extensively: when trying to prove some
property that is invariant by topological conjugacy (like the
existence of a free curve or a periodic point for $F$), it allows us
to consider just the case where the rotation set is the image of
$\rot(f)$ by some convenient element of $\GL(2,\Z)$.

\remark
  \label{rem:vertical}
  A particular case that will often appear is when $\rot(f)$ is a
  segment of rational slope. In that case, there exists a map
  $A\in\GL(2,\Z)$ such that $A\rot(f)$ is a vertical segment. Indeed,
  if $\rot(f)$ is a segment of slope $p/q$, then we can find $x,y\in
  \Z$ such that $px+qy=1$, and letting $$A=\left(\begin{matrix} p & -q \\
      y & x \end{matrix}\right)$$ it follows that $\det(A)=1$. Since
  $A(q,p)=(0,1)$, $A\rot(f)$ is vertical.
\remarkend

\subsection{The rotation set and periodic orbits.}

Recall that we say that a rational point $(p_1/q, p_2/q)\in \rot(f)$ is realized by a periodic orbit if there is $x\in \R^2$ such that $$f^q(x) = x + (p_1,p_2).$$
As we already mentioned in the introduction, a rational point in the
rotation set is not necessarily realized by a periodic orbit. However,
we have the following ``realization'' results (including the already
mentioned Theorem \ref{th:franks-area}). 

\begin{theorem}[Franks \cite{franks:reali-ext}]
  \label{th:reali-extremal}
  If a rational point of $\rot(f)$ is extremal, then it is realized by
  a periodic orbit.
\end{theorem}

\begin{theorem}[Franks \cite{franks:reali-int}]
  \label{th:reali-interior}
  If a rational point is in the interior of $\rot(f)$, then it is
  realized by a periodic orbit.
\end{theorem}

\begin{theorem}[Jonker \& Zhang \cite{jonker-zhang}]
  \label{th:jonker-zhang}
  If $\rot(f)$ is a segment with irrational slope and it contains a
  point of rational coordinates, then this point is realized by a
  periodic orbit.
\end{theorem}

We thank the referee for bringing the following generalization of the above theorem to our attention. It is stated for diffeomorphisms in \cite[p.\ 106]{lecalvez:asterisque}, but its proof remains valid for homeomorphisms using the results of \cite{lecalvez:translation} (see p.\ 9 of that article).

\begin{theorem}[Le Calvez]  If a rational point of $\rot(f)$ belongs to a line of irrational slope which bounds a closed half-plane that contains $\rot(f)$, then this point is realized by a periodic orbit.
\end{theorem}

\subsection{Curves and lines.}\label{sec:curves}

We denote by $I$ the interval $[0,1]$. A \emph{curve} on a manifold
$M$ is a continuous map $\gamma\colon I\to M$. As usual, we represent
by $\gamma$ both the map and its image, as it should be clear from the
context which is the case.

We say that the curve $\gamma$ is \emph{closed} if
$\gamma(0)=\gamma(1)$, and \emph{simple} if the restriction of the map
$\gamma$ to the interior of $I$ is injective. If $\gamma$ is a closed
curve, we say it is \emph{essential} if it is homotopically
non-trivial.

\definition
  A curve $\gamma\subset M$ is \emph{free} for $F$ if
  $F(\gamma)\cap\gamma=\emptyset$. We say that $F$ has the \emph{curve
    intersection property} if there are no free essential simple
  closed curves.
\definitionend
\remark
  For convenience, from now on by a \emph{free curve} for $F$ we will
  usually mean an essential simple closed curve that is free for $F$,
  unless otherwise stated.
\remarkend

By a \emph{line} we mean a proper topological embedding $\ell\colon \R
\to \R^2$. Again, we use $\ell$ to represent both the function and its
image. Note that lines are oriented.

\definition
  Given $(p,q)\in \ZCP$, a $(p,q)$-line in $\R^2$ is a line $\ell$
  such that there exists $\tau>0$ such that
  $\ell(t+\tau)=T_1^pT_2^q\ell(t)$ for all $t\in \R$, and such that
  its projection to $\T^2$ by $\pi$ is a simple closed curve. A
  $(p,q)$-curve in $\T^2$ is the projection by $\pi$ of a
  $(p,q)$-line. We will say that a simple closed curve is
  $\emph{vertical}$ if it is either a $(0,1)$-curve or a
  $(0,-1)$-curve. Similarly, a line will be called vertical if it is a
  $(0,1)$-line or a $(0,-1)$-line.
\definitionend

\remark
  If $\ell$ is a line in $\R^2$ that is invariant by $T_1^pT_2^q$,
  then $\pi(\ell)$ is always a closed curve in $\T^2$. We are
  requiring that this curve be simple to call $\ell$ a $(p,q)$-line.
  Conversely, if $\gamma$ is an essential simple closed curve in
  $\T^2$, taking a lift $\tilde{\gamma}\colon I\to\R^2$, we have that
  $(p,q) = \tilde{\gamma}(1)-\tilde{\gamma}(0)$ is an integer point
  independent of the choice of the lift. The curve $\tilde{\gamma}$
  can be extended naturally to $\R$ by
  $\tilde{\gamma}(t+n)=\tilde{\gamma}(t) + n(p,q)$, if $n\in \Z$ and
  $t\in [0,1]$; in this way we obtain a $(p,q)$-line that projects to
  $\gamma$. Since $\gamma$ is simple and essential, it is not hard to
  see that $p$ and $q$ must be coprime. With some abuse of notation,
  we will say that the $(p,q)$-line $\ell$ is a lift of $\pi(\ell)$.
\remarkend

\remark
  \label{rem:pq-bounded}
  Note that any $(p,q)$-line is contained in a strip bounded by two
  straight lines of slope $q/p$. In particular, if $\ell$ is a
  vertical line, there is $M>0$ such that $\pr_1(\ell)\subset [-M,M]$.
\remarkend

Given a line $\ell$ in $\R^2$, there are exactly two connected
components of $\R^2\setminus \ell$. Using the orientation of $\ell$,
we may define the \emph{left} and the \emph{right} components, which
we denote by $L\ell$ and $R\ell$. We also denote by $\ol L\ell$ and
$\ol R\ell$ their respective closures, which correspond to $L\ell\cup
\ell$ and $R\ell\cup \ell$.

Given two lines $\ell_1$ and $\ell_2$ in $\R^2$, we write $\ell_1<\ell_2$ if $\ell_1\subset L\ell_2$ and $\ell_2\subset
R\ell_1$. With an abuse of notation we will write $\ell_1\leq \ell_2$
when $\ell_1\subset \ol L\ell_2$, which means that the lines may
intersect but only ``from one side''. We will mostly use this relation to compare $(p,q)$-lines of the same type. Note that if $\ell_1$ and $\ell_2$ are disjoint $(p,q)$-lines, then either $\ell_1<\ell_2$ or
$\ell_2<\ell_1$ (but not both). In this context, we denote by $S(\ell_1,\ell_2)$
the strip $L\ell_2\cap R\ell_1$, and by $\ol S(\ell_1,\ell_2)$ its
closure.

\remark
  If $f\in\homeo(\R^2)$ preserves orientation, then $f$ preserves
  order: if $\ell_1<\ell_2$, then $f(\ell_1)<f(\ell_2)$.
\remarkend

\subsection{Brouwer lines.}

A Brouwer line for $h\in\homi(\R^2)$ is a line $\ell$ in $\R^2$ such
that $\ell< h(\ell)$ (sometimes Brouwer lines are not assumed to be oriented, but our lines are oriented). The classic Brouwer Translation Theorem
guarantees the existence of a Brouwer line through any point of $\R^2$
for any fixed-point free, orientation-preserving homeomorphism (see
\cite{brouwer:translation, kerekjarto:PB, fathi:translation,
  franks:translation, guillou:PB}; also see
\cite{lecalvez:translation} for a powerful equivariant version).  The
following will be much more useful for our purposes:

\begin{theorem}[Guillou \cite{guillou-preprint}]
  \label{th:guillou}
  Let $f$ be a lift of an orientation-preserving $F\in \homeo(\T^2)$,
  and suppose $f$ has no fixed points. Then $f$ has a Brouwer
  $(p,q)$-line, for some $(p,q)\in\ZCP$.
\end{theorem}

The following lemma is particularly useful when there is a Brouwer
$(p,q)$-line:
\begin{lemma}
  \label{lem:rot-brouwer}
  Let $S$ be a closed semiplane determined by a straight line
  containing the origin, and for $y\in\R^2$ denote by $S_y=\{w+y:w\in
  S\}$ its translate by $y$. Suppose that $x\in\R^2$ is such that for
  some $y$, $$f^n(x)\in S_y \text{ for all } n>0.$$ Then
  $\rot(f,x)\subset S$.  Moreover, if for all $x\in\R^2$ there is $y$
  such that the above holds, then $\rot(f)\subset S$.
\end{lemma}

\proc{Proof.}
  Let $\phi\colon \R^2\to\R$ be a linear functional such that
  $$S_y=\{w\in\R^2 : \phi(w)\geq \phi(y)\}.$$ Given $x$ such that
  $f^n(x)\in S_y$ for all $n>0$, we have then
  $$\phi\left(\frac{f^n(x)-x}{n}\right)=\frac{\phi(f^n(x))-\phi(x)}{n}\geq
  \frac{\phi(y) - \phi(x)}{n}\to 0.$$ Thus, if $z$ is the limit of a
  sequence of the form $(f^{n_i}(x)-x)/n_i$, then $\phi(z)\geq 0$.
  This implies that $\rot(f,x)\subset S$. The other claim follows from
  Theorem \ref{th:m-z-rot}.
\ep\medbreak

\remark
  \label{rem:semiplane}
  If $f$ has a Brouwer $(p,q)$-line $\ell$, the above lemma and Remark
  \ref{rem:pq-bounded} imply that $\rot(f)$ is contained in one of the
  closed semiplanes determined by the straight line $L$ of slope $q/p$
  through the origin. Indeed, by Remark \ref{rem:pq-bounded}, we can choose one of those semiplanes $S$ such that $S_y\supset R\ell\supset S_{y'}$ for some $y,y'\in \R^2$. Given $k\in \N$, we know that $\ell< f^k(\ell)$, so that $S_y \supset f^k(R\ell) \supset f^k(S_{y'})$. Given $x\in \R^2$, there is $(m,n)\in\Z^2$ such that $x+(m,n)\in S_{y'}$, and it follows  $f^k(x)\in S_{y-(m,n)}$ for all $k\in \N$. Thus by the above lemma, $\rho(f)\subset S$.
\remarkend

\subsection{The rotation set and free curves.}
Besides the existence of periodic orbits, other practical dynamical
information that can be obtained from the rotation set is the
existence of free curves. Recall that an interval with rational
endpoints $[p/q, p'/q']$ is a Farey interval if $qp' - pq' = 1$. The
following result was proved by Kwapisz for diffeomorphisms, and by
Beguin, Crovisier, LeRoux and Patou for homeomorphisms.

\begin{theorem}[\cite{kwapisz:degeneracy}, \cite{beg-et-al}]
  \label{th:beg-et-al}
  Suppose there exists a Farey interval $[p/q,p'/q']$ such that
  $$\pr_1\rot(f)\subset \left(\frac{p}{q},\frac{p'}{q'}\right).$$
  Then there exists a simple closed $(0,1)$-curve $\gamma$ in $\T^2$
  such that all the curves $\gamma$, $F(\gamma)$, $F^2(\gamma)$,
  $\ldots$, $F^{q+q'-1}(\gamma)$ are mutually disjoint.  In
  particular, if $\pr_1\rot(f)\cap \Z=\emptyset$, then $F$ has a free
  $(0,1)$-curve.
\end{theorem}

\subsection{The wedge.}
We now define an operation between $(p,q)$-lines, which is a fundamental
tool in the proof of Theorem A.

Recall that a Jordan domain is an open topological disk bounded by a simple
closed curve, and recall the following

\begin{theorem}[Ker\'ekj\'art\'o, \cite{kerekjarto:jordan}]
  If $U_1$ and $U_2$ are two Jordan domains in the two-sphere, then
  each connected component of $U_1\cap U_2$ is a Jordan domain.
\end{theorem}

Identifying $\R^2\cup\{\infty\}$ with the two-sphere
$\mathbb{S}^2$, we may regard lines in $\R^2$ as simple closed curves
in $\mathbb{S}^2$ containing $\infty$.  Recall that any two $(p,q)$-lines are oriented in the same way; that is, the intersection of their left sides is not contained in a strip. Thus if $\ell_1$ and $\ell_2$ are $(p,q)$-lines, it is easy to see that $L\ell_1\cap L\ell_2$ has a unique connected
component containing $\infty$ in its boundary. By the above theorem, this boundary is a simple closed curve, so that it corresponds to a
line in $\R^2$.  One can easily see that this new line is also a
$(p,q)$-line, if oriented properly. This motivates the following

\definition
  Given two $(p,q)$-lines $\ell_1$ and $\ell_2$ in $\R^2$, their
  \textit{wedge} $\ell_1\wedge\ell_2$ is the line defined as the
  boundary of the unique unbounded connected component of $L\ell_1\cap
  L\ell_2$, oriented so that this component corresponds to
  $L(\ell_1\wedge\ell_2)$.
\definitionend

\begin{figure}[!ht]
  \centering{\resizebox{0.6\textwidth}{!}{\includegraphics{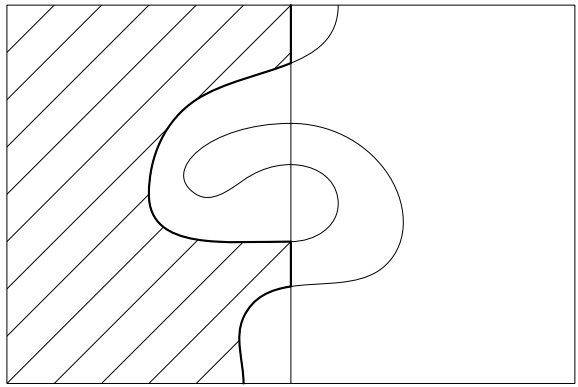}}}
  \caption{The wedge of two lines} \label{fig:wedge0}
\end{figure}

\remark
  This operation is called `join' in \cite{beg-et-al} and denoted by
  $\vee$.
\remarkend

We denote the wedge of multiple lines $\ell_1,\dots,\ell_n$ by
$$\ell_1\wedge\ell_2\wedge\cdots\wedge
\ell_n=\bigwedge_{i=1}^n\ell_i.$$ This is well defined because the
wedge is commutative and associative.
The following proposition resumes the interesting properties of the
wedge.

\begin{proposition}
  \label{prop:wedge}
  The wedge is commutative, associative and idempotent. Furthermore,
  \begin{enumerate}
  \item The wedge of $(p,q)$-lines is a $(p,q)$-line;
  \item If $h\in \homeo(\R^2)$ is a lift of a torus homeomorphism,
    then $h(\ell_1\wedge \ell_2)=h(\ell_1)\wedge h(\ell_2)$;
  \item $\ell_1\wedge \ell_2 \leq \ell_1$ and $\ell_1\wedge\ell_2 \leq
    \ell_2$;
  \item If $\ell_1<\ell_2$ and $\xi_1<\xi_2$, then $\ell_1\wedge\xi_1<
    \ell_2\wedge\xi_2$;
  \item The wedge of Brouwer lines is a Brouwer line.
  \end{enumerate}
\end{proposition}

\section{Realizing periodic orbits}

In this section we prove Theorem A. From now on we assume that
$F\in\homi(\T^2)$ and $f$ is a lift of $F$. The following result will
be essential in the proof.
\begin{proposition}
  \label{pro:CIP-iterado}
  Suppose $F^n$ has a free $(p,q)$-curve for some $n\geq 1$. Then $F$
  has a free $(p,q)$-curve.
\end{proposition}

We will also use the following lemmas, the proofs of which are
postponed to the end of this section.

\begin{lemma}
  \label{lem:CIP-1}
  Suppose $f^n$ has a Brouwer $(0,1)$-line, for some $n\in \N$. Then
  $f$ has a Brouwer $(0,1)$-line.
\end{lemma}

\begin{lemma}
  \label{lem:CIP-2}
  Suppose that some lift $f$ of $F$ has a Brouwer $(0,1)$-line. Then,
  either $F$ has a free $(0,1)$-curve, or $\max(\pr_1\rot(f))\geq 1$.
\end{lemma}

The next lemma is essentially Lemma 3 of \cite{bestvina-handel}; we
include it here for the sake of completeness.
\begin{lemma}
  \label{lem:CIP-3}
  Suppose $F$ has a free $(0,1)$-curve. Then for any lift $f$ of $F$
  there is $k\in \Z$ such that $\pr_1\rot(f)\subset [k,k+1]$.
\end{lemma}

\subsection{Proof of Proposition \ref{pro:CIP-iterado}.}

Conjugating the involved maps by an element of $\GL(2,\Z)$, we may
assume that $(p,q) = (0,1)$ (see Lemma \ref{lem:conjuga} and the
remark below it).

Suppose $F^n$ has a free $(0,1)$-curve for some $n\geq 2$. Then, any
lift to $\R^2$ of this curve, with a possible reversion of
orientation, is a vertical Brouwer line for $f^n$. Thus, Lemma
\ref{lem:CIP-1} implies that there is a vertical Brouwer line for $f$.
This holds for any lift $f$ of $F$.

If $\pr_1\rot(f)\cap \Z=\emptyset$, then by Theorem \ref{th:beg-et-al}
there is a free $(0,1)$-curve for $F$, and we are done.

Otherwise, let $k\in \pr_1\rot(f)\cap \Z$, and consider the lift
$f_0=T_1^{-k}f$ of $F$. By Remark \ref{pro:rot}, it is clear that
$0\in\pr_1\rot(f_0)$. On the other hand, as we already saw, $f_0$ has
a vertical Brouwer line $\ell$.

Assume $\ell$ is a Brouwer $(0,1)$-line. Then by Lemma
\ref{lem:CIP-2}, either $F$ has a free $(0,1)$-curve, or
$\max(\pr_1\rot(f_0))\geq 1$. In the latter case, since we know that
$0\in \rot(f_0)$, it follows from connectedness that
$\pr_1\rot(f_0)\supset [0,1]$. But this implies that
$$\pr_1\rot(f_0^n) \supset [0,n] \supset [0,2],$$ which contradicts
Lemma $\ref{lem:CIP-3}$ (since $F^n$ has a free $(0,1)$-curve).  Thus
the only possibility is that $F$ has a free $(0,1)$-curve.

If $\ell$ is a Brouwer $(0,-1)$-line, then using the previous argument
with $f_0^{-1}$ instead of $f_0$, we see that $F^{-1}$ (and thus $F$)
has a free $(0,-1)$-curve $\gamma$; and inverting the orientation of
$\gamma$ we get a free $(0,1)$-curve for $F$. This completes the
proof.  \ep

\subsection{Proof of Theorem A.}
Suppose $F\in \homi(\T^2)$ has the curve intersection property and
$\rot(f)$ has empty interior, where $f$ is a lift of $F$. We have
three cases.

\subsubsection{$\rot(f)$ is a single point.}
In this case, the unique point of $\rot(f)$ is extremal; and if it is
rational, Theorem \ref{th:reali-extremal} implies that it is realized
by a periodic orbit of $F$.

\subsubsection{$\rot(f)$ is a segment of irrational slope.}
In this case $\rot(f)$ contains at most one rational point and, by
Theorem \ref{th:jonker-zhang}, this point is realized by a periodic
orbit.

\remark
  \label{rem:jonker-zhang}
  Theorem $\ref{th:guillou}$ provides a simple way of proving this as
  well. In fact, it suffices to consider the case where the unique
  rational point in $\rot(f)$ is the origin, and to show that in this
  case $f$ has a fixed point. If the origin is an extremal point, this
  follows from Theorem \ref{th:reali-extremal}. If the origin is
  strictly inside the rotation set, then there is only one straight
  line through the origin such that $\rot(f)$ is contained in one of
  the closed semiplanes determined by the line. This unique line is
  the one with the same slope as $\rot(f)$, which is irrational. If
  $f$ has no fixed points, then Theorem $\ref{th:guillou}$ implies
  that $f$ has a Brouwer $(p,q)$-line for some $(p,q)\in \ZCP$; but
  then our previous claim contradicts Remark \ref{rem:semiplane}.
\remarkend

\subsubsection{$\rot(f)$ is a segment of rational slope.}
Fix a rational point $(p_1/q,p_2/q)\in \rot(f)$.  Recall that this
point is realized as the rotation vector of a periodic orbit of $F$ if
and only if $g = T_1^{-p_1}T_2^{-p_2}f^q$ has a fixed point. Note that
$(0,0)\in\rot(g)$, and $g$ is a lift of $F^q$. Moreover,
$$\rot(g) = T_1^{-p_1}T_2^{-p_2}(q\cdot\rot(f)),$$
which is a segment of rational slope containing the origin.
Conjugating $g$ by an element of $\GL(2,\Z)$, we may
assume that $\rot(g)$ is a vertical segment containing the origin  (\cf Remark \ref{rem:vertical}).

We will show by contradiction that $g$ has a fixed point. Suppose this
is not the case. Then, by Theorem $\ref{th:guillou}$, $g$ has a
Brouwer $(p,q)$-line $\ell$, for some $(p,q)\in \ZCP$. Moreover,
$(0,0)$ must be strictly inside $\rot(g)$ (\ie it cannot be extremal,
since otherwise $g$ would have a fixed point by Theorem
\ref{th:reali-extremal}), and by Remark \ref{rem:semiplane} this
implies that $\ell$ is a vertical Brouwer line.

Assume $\ell$ is a $(0,1)$-line (if it is a $(0,-1)$-line, we may
consider $g^{-1}$ instead of $g$ and use a similar argument). Since
$\diam(\pr_1 \rot(g)) = 0$, Lemma \ref{lem:CIP-2} implies that $F^q$,
the map lifted by $g$, has a free $(0,1)$-curve; but then by
Proposition \ref{pro:CIP-iterado}, $F$ has a free curve, contradicting
the curve intersection property. This concludes the proof.\ep

\subsection{Proof of Lemma \ref{lem:CIP-1}.}

Let $\ell$ be a Brouwer $(0,1)$-line for $f^n$, for some $n>1$. We
will show that there is a Brouwer $(0,1)$-line for $f^{n-1}$; by
induction, it follows that there is a Brouwer $(0,1)$-line for $f$.

We know that $\ell<f^n(\ell)$. Let $\xi$ be a $(0,1)$-line such that
$\ell<\xi < f^n(\ell)$, and define $$\ell' = \xi \wedge
\bigwedge_{i=1}^{n-1} f^i(\ell).$$

By Proposition \ref{prop:wedge}, $\ell'$ is still a $(0,1)$-line. We
claim that it is a Brouwer line for $f^{n-1}$. In fact,
\begin{align*}
  f^{n-1}(\ell') &= f^{n-1}(\xi)\wedge
  \bigwedge_{i=1}^{n-1}f^{n-1}(f^i(\ell))\\
  &= f^{n-1}(\xi)\wedge \bigwedge_{i=1}^{n-1}f^{i-1}(f^n(\ell))\\
  &= f^{n-1}(\xi)\wedge f^n(\ell) \wedge \bigwedge_{i=2}^{n-1}
  f^{i-1}(f^n(\ell))\\
  &= f^{n-1}(\xi)\wedge f^n(\ell) \wedge \bigwedge_{i=1}^{n-2}
  f^i(f^n(\ell)).
\end{align*}

Using the facts that $$f^{n-1}(\ell)<f^{n-1}(\xi),\,\quad\, 
\xi<f^n(\ell),\,\,\text{ and }\,\, f^i(\ell)<f^i(f^n(\ell)),$$ and Proposition
\ref{prop:wedge}, we see that $$\ell' = \xi\wedge \bigwedge_{i=1}^{n-1} f^i(\ell) = f^{n-1}(\ell)\wedge \xi \wedge \bigwedge_{i=1}^{n-2} f^i(\ell) < f^{n-1}(\ell'),$$ 
so that $\ell'$ is a Brouwer $(0,1)$-line for $f^{n-1}$. This concludes the proof.

\subsection{Proof of Lemma \ref{lem:CIP-2}.}
Let $\ell$ be a Brouwer $(0,1)$-line for $f$. We consider two cases.

\subsubsection{Case 1. For all $n>0$, $f^n(\ell)\nless T_1^n(\ell)$.}
In this case, for each $n>0$ we can choose $x_n\in\ell$ such that
$f^n(x_n)\in \ol R(T_1^n\ell)$. From the fact that $\ell$ is a
$(0,1)$-line we also know that $\pr_1(\ell)\subset [-M,M]$ for some
$M>0$, and therefore $$\pr_1(T_1^n(\ell))\subset [-M+n, M+n].$$ Hence,
$$\pr_1(f^n(x_n))\in \pr_1(R(T_1^n(\ell))) \subset [-M+n,\infty).$$
It then follows that $$\pr_1\left(\frac{f^n(x_n)-x_n}{n}\right)\geq
\frac{(-M+n) - M}{n} = -2\frac{M}{n} + 1 \xrightarrow[]{n\to\infty}
1,$$ and by the definition of rotation set this implies that some
point $(x,y)\in \rot(f)$ satisfies $x\geq 1$; \ie
$\max\pr_1(\rot(f))\geq 1$.

\subsubsection{Case 2. $f^n(\ell)<T_1^n\ell$ for some $n>0$.}

We will show that in this case $F$ has a free $(0,1)$-curve. The idea
is similar to the proof of Lemma \ref{lem:CIP-1}. Let $n$ be the
smallest positive integer such that $f^n(\ell)<T_1^n\ell$. If $n=1$,
we are done, since $\ell<f(\ell)<T_1\ell$, so that $\ell$ projects to
a free $(0,1)$-curve for $F$.

Now assume $n>1$. We will show how to construct a new Brouwer
$(0,1)$-line $\beta$ for $f$ such that
$f^{n-1}(\beta)<T_1^{n-1}\beta$. Repeating this argument $n-1$ times,
we end up with a $(0,1)$-line $\ell'$ such that
$\ell'<f(\ell')<T_1(\ell')$, so that $\ell'$ projects to a free curve,
completing the proof.

Let $\xi$ be a $(0,1)$-curve such that
\begin{equation}
  \label{eq:dem-wedge1}
  f^n(\ell)<\xi<T_1^n\ell
\end{equation}
We may choose $\xi$ such that $\xi < f(\xi)$, by taking it close enough
to $f^n(\ell)$. This is possible because $f^n(\ell)< f(f^n(\ell))$, and
these two curves are separated by a positive distance, since they both project to
$(0,1)$-curves in $\T^2$. Thus $\xi$ is also a Brouwer $(0,1)$-line
for $f$.

Define $$\beta= \xi \wedge \bigwedge_{i=1}^{n-1} T_1^{n-i}f^i(\ell).$$
Let us see that $f^{n-1}(\beta)<T_1^{n-1}(\beta)$. Since $T_1$ commutes
with $f$, we have
\begin{align*}
  f^{n-1}(\beta) &= f^{n-1}(\xi)\wedge \bigwedge_{i=1}^{n-1}
  f^{n-1}T_1^{n-i}f^i(\ell)\\
  &= f^{n-1}(\xi)\wedge \bigwedge_{i=1}^{n-1} f^{i-1}T_1^{n-i}f^n(\ell)\\
  &= f^{n-1}(\xi)\wedge T_1^{n-1}f^n(\ell) \wedge
  \bigwedge_{i=2}^{n-1}f^{i-1}T_1^{n-i}f^n(\ell)
\end{align*}
By (\ref{eq:dem-wedge1}), we also have
\begin{itemize}
\item $f^{n-1}(\xi)<f^{n-1}(T_1^n\ell),$
\item $T_1^{n-1}f^n(\ell)< T_1^{n-1}\xi,$ and
\item $f^{i-1}T_1^{n-i}f^n(\ell)< f^{i-1}T_1^{n-i}T_1^n\ell$.
\end{itemize}
Using these facts and Proposition \ref{prop:wedge} we see that
\begin{align*}
  f^{n-1}(\beta) &< f^{n-1}(T_1^n\ell)\wedge T_1^{n-1}\xi \wedge
  \bigwedge_{i=2}^{n-1}f^{i-1}T_1^{n-i}T_1^n(\ell)\\
  &= T_1^{n-1}(T_1f^{n-1}(\ell)) \wedge
  T_1^{n-1}(\xi)\bigwedge_{i=2}^{n-1}T_1^{n-1}(T_1^{n-(i-1)}f^{i-1}(\ell))\\
  &= T_1^{n-1}\left(\xi\wedge T_1f^{n-1}(\ell) \wedge
    \bigwedge_{i=1}^{n-2}T_1^{n-i}f^i(\ell)\right)\\
  &= T_1^{n-1}\left(\xi\wedge\bigwedge_{i=1}^{n-1}
    T_1^{n-i}f^i(\ell)\right)\\
  &= T_1^{n-1}(\beta).
\end{align*}

Since $\beta$ is a Brouwer $(0,1)$-line, this completes the proof.\ep

\subsection{Proof of Lemma \ref{lem:CIP-3}.}
Let $f$ be a lift of $F$, and suppose that $F$ has a free $(0,1)$-curve. A lift of this curve is a $(0,1)$-line $\ell$ such that $T_1^nf(\ell)\cap \ell=\emptyset$ for any $n\in\N$. In particular, given $n\in \N$, either $\ell < T_1^nf(\ell)$ or $T_1^nf(\ell)<\ell$ (and in the latter case reversing the orientation of $\ell$ we get a Brouwer $(0,-1)$-line for $T_1^nf$). In either case, Remark \ref{rem:semiplane} implies that $\rot(T_1^nf)$ is contained in one of the semiplanes $\{(x,y):x\geq 0\}$ or $\{(x,y):x\leq 0\}$, and by Remark \ref{pro:rot} this means that $\rot(f)$ is contained in one of the semiplanes $\{(x,y):x\geq -n\}$ or $\{(x,y):x\leq -n\}$. Thus $n$ cannot be an interior point of $\pr_1 \rot(f)$. Since this holds for all $n$, and $\pr_1\rot(f)$ is an interval, there is $k$ such that $\pr_1\rot(f)\subset[k,k+1]$. \ep

\section{Free curves and fixed points}
\label{ch:theoremB}

In this section we prove Theorem B. For the proof, we have two main
cases. The first one is when the rotation set is either a segment of
rational slope or a single point; this is dealt with Theorem A and the
results stated in \S\ref{sec:prelim}. The second case is when the
rotation set is a segment of irrational slope. In that case, the main
idea is to find $A\in \GL(2,\R)$ such that $A\rot(f)$ has no integers
in the first or second coordinate, so that we may apply directly
Theorem \ref{th:beg-et-al} to $AfA^{-1}$ (\cf Lemma
\ref{lem:conjuga}). In fact, using this argument we obtain the
following more general result:

\begin{theorem}
  \label{th:irrational}
  Suppose $\rot(f)$ is a segment of irrational slope with no rational
  points. Then for each $n>0$ there is an essential simple closed
  curve $\gamma$ such that $\gamma,\,F(\gamma),\dots,\,F^n(\gamma)$
  are pairwise disjoint.
\end{theorem}

\remark No example is known of a homeomorphism that meets the hypothesis of the above theorem. In fact, if a conjecture of Franks and Misiurewicz turns out to be true (see \cite{f-m}), then no such example exists. This theorem could be a (very small) step towards this conjecture.
\remarkend

The problem of finding the map $A$ previously mentioned is mainly
arithmetic, and we consider it first. In the next section, we briefly
recall a few facts about continued fractions that will be needed in
the proof; in \S\ref{sec:arithmetic} we prove the two arithmetic
lemmas that allow us to find the map $A$; in \S\ref{sec:Bproof1} we
prove Theorem \ref{th:irrational}; finally, in \S\ref{sec:Bproof2} we
complete the proof of Theorem B.

\subsection{Continued fractions.}

Given an integer $a_0$ and positive integers $a_1,\dots,a_n$, we
define
$$[a_0; a_1,\dots,a_n]=a_0+ \cfrac{1}{a_1+ \cfrac{1}{\ddots
    \genfrac{}{}{0pt}{}{}{{}+\,\cfrac{1}{a_n}}}}.$$

Given $\alpha\in\R$, define $\{\alpha_n\}$ and $\{a_n\}$ recursively
by $a_0=\floor{\alpha}$, $\alpha_0 = \alpha-a_0$, and
$$a_{n+1} = \floor{\alpha_n^{-1}}, \,\, \alpha_{n+1} = \alpha_n^{-1}-a_{n+1},$$
whenever $\alpha_n\neq 0$. This gives the continued fractions
representation of $\alpha$: If $\alpha$ is rational, we get a finite
sequence $a_0,\dots,a_n$, and $$\alpha=[a_0;a_1,\dots,a_n].$$ If
$\alpha$ is irrational, then the sequence is infinite and
$$\alpha=[a_0;a_1,a_2,\dots]\doteq
\lim_{n\to\infty}[a_0;a_1,\dots,a_n].$$

The rational number $p_n/q_n=[a_0;a_1\dots,a_n]$ is called the $n$-th
\emph{convergent} to $\alpha$. Convergents may be regarded as the
``best rational approximations'' to $\alpha$, in view of the following
properties (see, for instance, \cite{har-wri})
\begin{proposition}
  \label{pro:cfrac}
  If $p_n/q_n$ are the convergents to $\alpha$, then
  \begin{enumerate}
  \item $\{q_n\}$ is an increasing sequence of positive integers, and
    $$\frac{1}{q_n+q_{n+1}} < (-1)^n(\alpha q_n - p_n) <
    \frac{1}{q_{n+1}}.$$
  \item $\frac{p_{2n}}{q_{2n}} < \frac{p_{2n+2}}{q_{2n+2}} < \alpha <
    \frac{p_{2n+3}}{q_{2n+3}}<\frac{p_{2n+1}}{q_{2n+1}}$.
  \item $p_{n+1}q_n - p_nq_{n+1}=(-1)^n$
  \end{enumerate}
\end{proposition}

A Farey interval is a closed interval with rational endpoints
$[p/q,p'/q']$ such that $p'q-q'p=1$. Note that two consecutive
convergents give Farey intervals.

\begin{proposition}
  \label{prop:farey}
  Let $[p/q,p'/q']$ be a Farey interval. Then
  $$\left[\frac{p+p'}{q+q'},\frac{p'}{q'}\right]\, \text{ and }
  \left[\frac{p}{q},\frac{p+p'}{q+q'}\right]$$ are Farey intervals,
  and if $p''/q''\in (p/q,p'/q')$, then $q''\geq q + q'$.
\end{proposition}

\subsection{Arithmetic lemmas.}
\label{sec:arithmetic}

Define the vertical and horizontal inverse \textit{Dehn twists} by
$$D_1\colon (x,y)\mapsto (x-y,y),\,\,D_2\colon (x,y)\mapsto (x,y-x).$$
Let $Q$ be the set of vectors of $\R^2$ with positive coordinates; for
$u=(x,y)\in Q$, we denote by $\slo(u) = y/x$ the \emph{slope} of $u$.
Let $Q_1$ and $Q_2$ be the sets of elements of $Q$ having slope
smaller than one and greater than one, respectively.

\remark
  \label{rem:twist}
  Note the following simple properties
  \begin{enumerate}
  \item for $i=1,2$, $D_iQ_i=Q$; and if $u\in Q$, then $D_i^{-k}u\in
    Q_i$ for all $k>0$;
  \item for $i=1,2$, $\norm{D_iu}<\norm{u}$ if $u\in Q_i$;
  \item $\slo(D_2^ku) = \slo(u)-k$ and $\slo(D_1^ku)^{-1} =
    \slo(u)^{-1}-k$.
  \end{enumerate}
\remarkend

\begin{lemma}
  \label{lem:diag}
  Let $u$ and $v$ be elements of $Q$ with different slopes. Then there
  is $A\in \GL(2,\Z)$ such that
  \begin{enumerate}
  \item $\norm{Au}\leq \norm{u}$;
  \item $\norm{Av}\leq \norm{v}$;
  \item Both $Au$ and $Av$ are in $Q$, and either one of these points
    is on the diagonal, or one is in $Q_1$ and
    the other in $Q_2$.
  \end{enumerate}
\end{lemma}

\proc{Proof.}
  We first note that it suffices to consider the case where both $u$
  and $v$ are in $Q_1$.  Indeed, if one of the vectors is in $Q_i$ and
  the other is not (for $i=1$ or $2$), we can set $A=Id$; and if
  both $u$ and $v$ are in $Q_2$ then we may use $Su$ and $Sv$ instead,
  where $S$ is the isometry $(x,y)\mapsto (y,x)$.

\begin{figure}[!ht]
  \centering
  \includegraphics{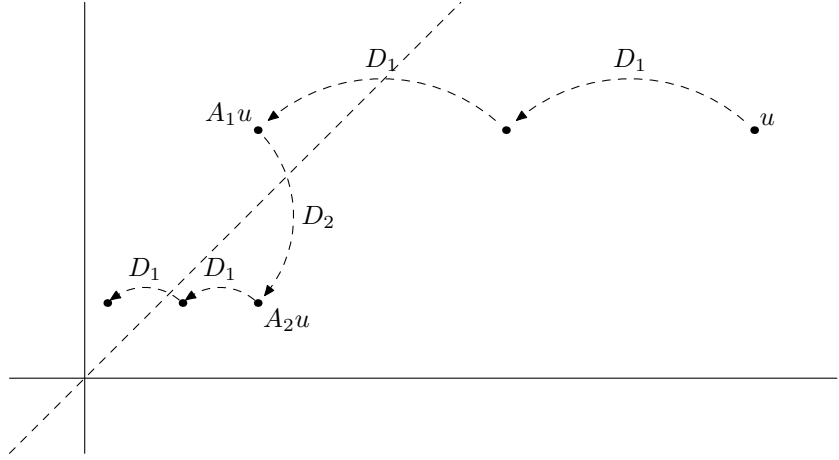}
  \caption{The sequence $A_iu$} \label{fig:twist0}
\end{figure}

  Given $u\in Q_1$, we define a sequence of matrices $A_n\in
  \SL(2,\Z)$ and integers $a_n$ by $A_0=I$, $a_0=0$, and recursively
  (see Figure \ref{fig:twist0})
  \begin{itemize}
  \item If $\slo(A_nu)=1$ stop the construction.
  \item $a_{n+1}$ is the smallest integer such that
    $D_i^{a_{n+1}}A_nu\notin Q_i$, where $i=2$ if $n$ is odd, $1$ if
    $n$ is even;
  \item $A_{n+1}=D_i^{a_{n+1}}A_n$.
  \end{itemize}

  In this way we get either an infinite sequence, or a finite sequence
  $A_1,\dots A_N$ such that $A_Nu$ lies on the diagonal and has
  positive coordinates. Furthermore, given $0\leq n<N$ if the sequence
  is finite, or $n\geq 0$ if it is infinite, we have
  \begin{enumerate}
  \item $A_nu\in Q_i$ where $i=2$ if $n$ is odd, $1$ if $n$ is even;
  \item If $\alpha_n=\slo(A_nu)^{(-1)^n}$, then $\alpha_n =
    \alpha_{n-1}^{-1} - a_n$
  \item $\norm{A_0u},\norm{A_1u},\dots$ is a decreasing sequence;
  \end{enumerate}

The first property is a consequence of the definition. The second
follows from Remark \ref{rem:twist}, since for odd $n$, we have
$$\alpha_n = \slo(D_1^{a_n}A_{n-1}u)^{-1} = \slo(A_{n-1}u)^{-1}-a_n =
\alpha_{n-1}^{-1} - a_n$$ while for even $n$,
$$\alpha_n = \slo(D_2^{a_n}A_{n-1}u) = \slo(A_{n-1}u)-a_n =
\alpha_{n-1}^{-1}-a_n.$$ The last property also follows from the
construction: if $A_nu=D_i^{a_n}A_{n-1}u$, then $D_i^kA_{n-1}u\in Q_i$
for all $0\leq k<a_n$, so that Remark \ref{rem:twist} implies that
$\norm{A_nu}<\norm{A_{n-1}u}$.

If $n$ is odd, $a_n$ is the smallest integer such that
$D_1^{a_n}A_{n-1}u\notin Q_1$, or equivalently (assuming that $n<N$ if
the sequence $a_n$ is finite), the smallest integer such that
$$\slo(D_1^{a_n}A_{n-1}u)=(\slo(A_{n-1}u)^{-1}-a_n)^{-1} =
(\alpha_{n-1}^{-1}-a_n)^{-1}\geq 1.$$ Since $A_{n-1}u\in Q_1$,
$\slo(A_{n-1}u)<1$ so that $\alpha_{n-1}^{-1}>1$, and $a_n>0$. Note
that $\alpha_n^{-1}$ cannot be an integer, since otherwise
$\slo(A_nu)=1$, which contradicts the fact that $A_nu\in Q_2$; thus
$$a_n=\floor{\alpha_{n-1}^{-1}}.$$
If $n$ is even, the above equation holds by a similar argument. One
easily sees from these facts that $\{\alpha_n\}$ coincides with the
sequence obtained in the definition of the continued fractions
expression of $\alpha_0=\slo(u)$, and thus $\{a_n\}$ coincides with
the continued fractions coefficients of $\slo(u)$.

Now given $v\in Q_1$ with $\slo(v)\neq\slo(u)$, define in the same way
as above sequences of positive integers $b_n$ and matrices
$B_n\in\SL(2,\Z)$ such that $B_0=I$, $b_n$ is given by the continued
fractions expression of $\slo(v)$, and $B_{n+1}=D_i^{b_{n+1}}B_n$
where $i=2$ if $n$ is odd, $1$ if $n$ is even. As before, we have that
$\norm{B_nv}$ is a (finite or infinite) decreasing sequence, and if it
is finite of length $N$ then $\slo(B_Nv)=1$. Also $B_nv\in Q_i$ where
$i=2$ if $n$ is odd, $1$ if $n$ is even (given that $n<N$ if the
sequence is finite).

\begin{figure}[!ht]
  \centering
  \includegraphics{figures/twisting1.1}
  \caption{Example}
  \label{fig:twist1}
\end{figure}

Since $u$ and $v$ have different slopes, the continued fractions
expressions of these slopes cannot coincide.  Thus there exists $m\geq 0$ such that
$a_0=b_0,\dots,a_m=b_m$, and either $b_{m+1}\neq a_{m+1}$, or exactly one of $a_{m+1}$ or $b_{m+1}$ is not defined. In the latter case, we may assume that $a_{m+1}$ is undefined (by swapping $u$ and $v$ if necessary) and this means that $\slo(A_mu)=1$ (from the previous construction), so that $A=A_m=B_m$ meets the required conditions. In the former case, we may further assume that $a_{m+1}<b_{m+1}$ (again, by swapping $u$ and $v$ if necessary).
This means that $A_k=B_k$ for $0\leq k\leq m$; so that if $m$ is even,
$A_mu\in Q_1$ and $A_mv\in Q_1$, but since $a_{m+1}<b_{m+1}$, it holds
that
$$A_{m+1}u = D_1^{a_{m+1}}A_mu\notin Q_1\,\text{ but }\, A_{m+1}v =
D_1^{a_{m+1}}B_mv\in Q_1;$$ that is, $D_1^{a_{m+1}}$ ``pushes'' $A_mu$
out of $Q_1$, while leaving $A_mv$ in $Q_1$ (see Figure
\ref{fig:twist1}).

Moreover, since $a_{m+1}$ is minimal with that property, either
$\slo(A_{m+1}u)=1$ or $A_{m+1}u\in Q_2$; and in either case $A_{m+1}u$
has positive coordinates.  By construction, it also holds that
$\norm{A_{m+1}v}\leq\norm{v}$ and $\norm{A_{m+1}u}\leq\norm{u}$.

If $m$ is odd, a similar argument holds, and we see that $A_{m+1}u\in
Q_2$ while $A_{m+1}v$ is either on the diagonal or in $Q_1$.

Starting from $u$ and $v$ in $Q_1$ we obtained $A=A_{m+1}\in\GL(2,\Z)$
such that $Au$ and $Av$ have positive coordinates and either one of
them is on the diagonal or one is in $Q_1$ and the other in $Q_2$;
furthermore, $\norm{Au}\leq \norm{u}$ and $\norm{Av}\leq \norm{v}$.
This completes the proof.
\ep\medbreak

\begin{lemma}
  \label{lem:size}
  Let $w=(x,y)$ be a vector with irrational slope. Then, for any
  $\epsilon>0$ there exists $A\in \SL(2,\Z)$ such that
  $\norm{Aw}<\epsilon$.
\end{lemma}

\proc{Proof.}
  Let $p_i/q_i$ be the convergents to $y/x$, and define $$ A_k =
  \left(
    \begin{matrix}
      (-1)^kp_k & (-1)^{k+1}q_k \\ p_{k+1} & -q_{k+1}
    \end{matrix}
  \right).$$ By Proposition \ref{pro:cfrac} we have that $\det A=1$,
  the sequence $q_1,q_2,\dots$ is increasing, and $$\abs{p_i -
    \frac{y}{x} q_i}< \frac{1}{q_{i+1}}$$ for all $i\geq 1$. Hence,
  $$\abs{\pr_1 A_kw}= \abs{x\left(p_k - \frac{y}{x} q_k\right)} <
  \frac{\abs{x}}{q_{k+1}},$$ and similarly $$\abs{\pr_2 A_kw} =
  \abs{x\left(p_{k+1}-\frac{y}{x} q_{k+1}\right)}<
  \frac{\abs{x}}{q_{k+2}};$$ Choosing $k$ large enough so that
  $q_{k+1} > \sqrt{2}\abs{x}\epsilon^{-1}$, we have
  $\norm{A_kw}<\epsilon$.
\ep\medbreak

\subsection{Proof of Theorem \ref{th:irrational}.}
\label{sec:Bproof1}
We divide the proof into two cases:

\subsubsection{The case $n=1$}. We first assume $n=1$, \ie we prove that $F$
has a free curve, assuming that $\rot(f)$ is a segment of irrational
slope containing no rational points. The problem is reduced, by Lemma
\ref{lem:conjuga} and Theorem \ref{th:beg-et-al}, to finding $A\in
\GL(2,\Z)$ such that the projection of $A\rot(f)$ to the first or the
second coordinate contains no integers.  Note that Lemma
\ref{lem:size} allows us to assume that $$\diam(\rot(f))<\epsilon <
\frac{1}{2\sqrt{5}}.$$ We may also assume that there are $m_1\in
\pr_1(\rot(f))\cap \Z$ and $m_2\in \pr_2(\rot(f))\cap\Z$, for
otherwise there is nothing to do.

Then using $T_1^{-m_1}T_2^{-m_2}f$ (which also lifts $F$) instead of
$f$, we have that the extremal points of $\rot(f)$ are in opposite
quadrants. By conjugating $f$ with a rotation by $\pi/2$, we may
assume that $\rot(f)$ is the segment joining $u=(-u_1,-u_2)$ and
$v=(v_1,v_2)$ where $v_i\geq 0$ and $u_i\geq 0$, $i=1,2$. From this,
and the fact that $\diam(\rot(f))<\epsilon$, it follows that
$$\norm{u}< \epsilon,\,\text{ and }\, \norm{v}<\epsilon.$$

\subsubsection*{Case 1. One of the points has a zero coordinate.}
\begin{figure}[ht]
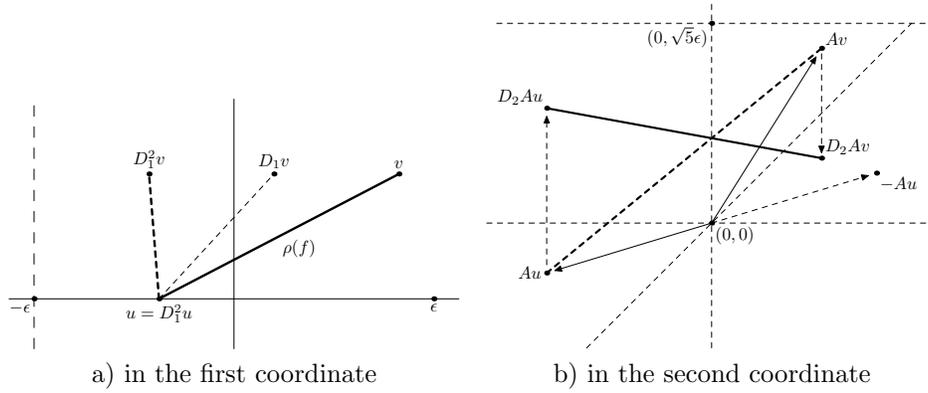

  \centering $\begin{array}{c@{\hspace{10pt}}c}
   \resizebox{6cm}{!}{\includegraphics{figures/twisting2.1}} &
    \resizebox{6cm}{!}{\includegraphics{figures/twisting3.1}} \\ 
    \mbox{a) in the first coordinate} & \mbox{b) in
      the second coordinate}
  \end{array}$
  \caption{Avoiding integers}
  \label{fig:twisting}
\end{figure}

It is clear that neither $u$ nor $v$ can have both coordinates equal
to $0$. Conjugating by an appropriate isometry in $\GL(2,\Z)$, we may
assume the generic case that $u=(-u_1,0)$, with $u_1>0$. Then $v_2>0$:
in fact if $v_2=0$ then $\rot(f)$ contains the origin, which is not
possible. Let $k>0$ be the greatest integer such that
$\pr_1D_1^kv>-u_1$, \ie $$k=\floor{\frac{u_1+v_1}{v_2}}.$$ Note that
$D_1^ku=u$, so that $D_1^k\rot(f)$ is the segment joining $u$ to
$D_1^kv$ (see Figure \ref{fig:twisting}a).
Moreover, $D_1^{k+1}v = (v_1',v_2)$ where $v_1' = v_1-(k+1)v_2 <
-u_1$. Thus $$\max \pr_1(D_1^{k+1}\rot(f)) = -u_1<0.$$ On the other
hand, $$\min \pr_1(D_1^{k+1}\rot(f)) \geq -u_1-v_1 > -2\epsilon>-1,$$
so that taking $A=D_1^{k+1}$ we have $\pr_1(\rot(AfA^{-1}))\subset
(-1,0)$. By Theorem \ref{th:beg-et-al}, it follows that $\tilde A F
\tilde A^{-1}$ has a free $(0,1)$ curve. Thus $F$ has a free curve.

\subsubsection*{Case 2. None of the points has a zero coordinate.}

In this case, $u\in -Q$ and $v\in Q$. Since the segment joining $u$ to
$v$ cannot contain the origin, $-u$ and $v$ are elements of $Q$ with
different slopes; thus Lemma \ref{lem:diag} implies that there is
$A\in \GL(2,\Z)$ such that $\norm{Au}\leq \norm{u}$, $\norm{Av}\leq
\norm{v}$, both $-Au$ and $Av$ are in $Q$, and either one of them lies
on the diagonal, or they are in opposite sides of the diagonal. This
means that $Au$ and $Av$ are both contained in one of the closed
semiplanes determined by the diagonal. By using $-A$ instead of $A$ if
necessary, we may assume that both are in the closed semiplane above
the diagonal, which is mapped by $D_2$ to the upper semiplane
$H=\{(x,y):y\geq 0\}$.  Note that $$\norm{D_2Au}\leq
\norm{D_2}\norm{Au}\leq \sqrt{5}\norm{u},$$ and similarly
$\norm{D_2Av}\leq \sqrt{5}\norm{v}$.  If $\pr_2(D_2Av)>0$ and
$\pr_2(D_2Au)>0$, then (see Figure \ref{fig:twisting}b) we have that
$$\pr_2\rot(f)\subset (0,\sqrt{5}\epsilon)\subset (0,1),$$ and by
Theorem \ref{th:beg-et-al} (as in Case 1) it follows that $F$ has a
free curve. On the other hand, if either of $D_2Av$ or $D_2Au$ has
zero second coordinate, the argument in Case 1 implies that $F$ has a
free curve.

This completes the proof when $\rot(f)$ has irrational slope and
$n=1$.

\subsubsection{The case $n>1$.} Note that when $\rot(f)$ has irrational
slope, $\rot(f^n)=n\rot(f)$ has irrational slope for all $n$.

Let $N=n!$. As we saw in the previous case, conjugating our maps by
some $A\in \GL(2,\Z)$, we may assume $\pr_1\rot(f^N)\cap\Z=\emptyset$;
thus $$\pr_1\rot(f^N)\subset (K,K+1)\,\text{ for some } K\in \Z,$$ and
by Remark \ref{pro:rot}, $$\pr_1\rot(f)\subset \left(\frac{K}{N},
  \frac{K+1}{N}\right).$$ Let $[p/q,p'/q']$ be the smallest Farey
interval containing $\pr_1(\rot(f))$.  We claim that $q+q'> n$. In
fact, if $q+q'\leq n$, then $[K/N,(K+1)/N]$ must be contained in one
of the smaller Farey intervals (see Proposition \ref{prop:farey})
$$\left[\frac{p}{q},\frac{p+p'}{q+q'}\right]\,\text{ or
}\,\left[\frac{p+p'}{q+q'},\frac{p'}{q'}\right].$$ This is because
$N(p+p')/(q+q')$ is an integer, so that it cannot be in the interior
of $[K, K+1]$. But we chose our Farey interval to be the smallest, so
$[p/q,p'/q']$ must as well be contained in one of these two intervals,
which is a contradiction. Thus $q+q'> n$, and Theorem
\ref{th:beg-et-al} guarantees that there is an essential simple closed
curve $\gamma$ such that its first $n$ iterates by $F$ are pairwise
disjoint. This completes the proof.

\subsection{Proof of Theorem B.}
\label{sec:Bproof2}

Assume that $\rot(f)$ has empty interior. We will show that either $F$
has a fixed point, or it has a free curve. There are several cases:
\begin{itemize}
\item $\rot(f)$ is a segment of irrational slope which contains no
  rational points. Then there is a free curve, by Theorem
  \ref{th:irrational}.
\item $\rot(f)$ is a segment of irrational slope containing a rational
  non-integer point $(p_1/q, p_2/q)$. By Lemma $\ref{lem:size}$ there
  exists $A\in \GL(2,\Z)$ such that $\diam(A\rot(f))<1/q$. One of the
  two coordinates of $A(p_1/q,p_2/q)$ must be non-integer. We assume
  $p'/q'=\pr_1A(p_1/q,p_2/q)\notin \Z$ (otherwise, we can conjugate
  $f$ with a rotation by $\pi/2$, as usual). Since $A(p_1,p_2)$ is an
  integer point, it follows that $q'\leq q$ (if we assume $p'/q'$ is
  irreducible). Thus,
  $$\pr_1\left(\rot(AfA^{-1})\right)=\pr_1(A\rot(f)) \subset \pr_1
  \left(\frac{p'}{q'}-\frac{1}{q},\frac{p'}{q'}+\frac{1}{q}\right).$$
  It is clear that the interval above contains no integers, so that
  $\tilde AF\tilde A^{-1}$ (and, consequently, $F$) has a free curve
  by Theorem $\ref{th:beg-et-al}$.
\item $\rot(f)$ is a segment of irrational slope with an integer
  point. Then $F$ has a fixed point by Theorem \ref{th:jonker-zhang}
  (see also Remark \ref{rem:jonker-zhang}).
\item $\rot(f)$ is a single point. Then either this point is integer,
  and $F$ has a fixed point by Theorem \ref{th:reali-extremal} or it
  is not integer, and $F$ has a free curve by Theorem
  \ref{th:beg-et-al}.
\item $\rot(f)$ is a segment of rational slope. Conjugating all the
  maps by an element of $\GL(2,\Z)$ we may assume it is a vertical
  segment; and with this assumption, if both $\pr_1(\rot(f))$ and
  $\pr_2(\rot(f))$ contain an integer, it follows that $\rot(f)$
  contains an integer point, and by Theorem A, $F$ has either a fixed
  point or a free curve. On the other hand, if either of the two
  projections contains no integer, Theorem \ref{th:beg-et-al} implies
  the existence of a free curve for $F$ as before.
\end{itemize}
This concludes the proof.

\acks
The results presented in this work are part of the doctoral thesis of
the second author, under the supervision of Enrique Pujals. We are
thankful to him for the constant encouragement and valuable
conversations. We are indebted to Patrice Le Calvez for insightful
discussions. We would also like to thank Sebasti\~ao Firmo and Meysam
Nassiri for many corrections and suggestions.

\bibliographystyle{amsalpha} 
\bibliography{tesis}
\end{document}